\documentclass[12pt]{article}
\usepackage[dvips]{epsfig}
\usepackage{amssymb}
\usepackage{amsthm}
\usepackage{amsmath}
\usepackage[T1]{fontenc}
\usepackage[utf8]{inputenc}
\usepackage{bm}
\usepackage{dcolumn}

\catcode`\@=11
\def\marginnote#1{}
\newcount\hour
\newcount\minute
\newtoks\amorpm
\hour=\time\divide\hour by60
\minute=\time{\multiply\hour by60 \global\advance\minute
by-\hour}\edef\standardtime{{\ifnum\hour<12
\global\amorpm={am}%
        \else\global\amorpm={pm}\advance\hour by-12 \fi
        \ifnum\hour=0 \hour=12 \fi
        \number\hour:\ifnum\minute<10
0\fi\number\minute\the\amorpm}}
\edef\militarytime{\number\hour:\ifnum\minute<10
0\fi\number\minute}

\def\draftlabel#1{{\@bsphack\if@filesw {\let\thepage\relax
   \xdef\@gtempa{\write\@auxout{\string
      \newlabel{#1}{{\@currentlabel}{\thepage}}}}}\@gtempa
   \if@nobreak \ifvmode\nobreak\fi\fi\fi\@esphack}
        \gdef\@eqnlabel{#1}}
\def\@eqnlabel{}
\def\@vacuum{}
\def\draftmarginnote#1{\marginpar{\raggedright\scriptsize\tt#1}}
\def\draft{\oddsidemargin -.5truein
        \def\@oddfoot{\sl preliminary draft \hfil
        \rm\thepage\hfil\sl\today\quad\militarytime}
        \let\@evenfoot\@oddfoot \overfullrule 3pt
        \let\label=\draftlabel
        \let\marginnote=\draftmarginnote

\def\@eqnnum{(\theequation)\rlap{\kern\marginparsep\tt\@eqnlabel}%
\global\let\@eqnlabel\@vacuum}  }

\def\numberbysection{\@addtoreset{equation}{section}
        \def\theequation{\thesection.\arabic{equation}}}

\def\underline#1{\relax\ifmmode\@@underline#1\else
 $\@@underline{\hbox{#1}}$\relax\fi}

\catcode`@=12
\relax

\numberbysection

\advance \topmargin by -\headheight
\advance \topmargin by -\headsep

\evensidemargin \oddsidemargin
\def\br{\begin{eqnarray}}
\def\er{\end{eqnarray}}
\def\be{\begin{equation}}
\def\ee{\end{equation}}

\def\({\left(}
\def\){\right)}

\relax

\def\tp0{\Theta_{+}^{(0)}}
\def\tm0{\Theta_{-}^{(0)}}

\def\f#1#2#3 {f^{#1#2}_{#3}}

\def\win1{{\sf w_{1+\infty}}}

\def\Win1{{\sf W_{1+\infty}}}

\def\rlx{\relax\leavevmode}
\def\inbar{\vrule height1.5ex width.4pt depth0pt}
\def\IZ{\rlx\hbox{\sf Z\kern-.4em Z}}
\def\IR{\rlx\hbox{\rm I\kern-.18em R}}
\def\IC{\rlx\hbox{\,$\inbar\kern-.3em{\rm C}$}}
\def\IN{\rlx\hbox{\rm I\kern-.18em N}}
\def\IO{\rlx\hbox{\,$\inbar\kern-.3em{\rm O}$}}
\def\IP{\rlx\hbox{\rm I\kern-.18em P}}
\def\IQ{\rlx\hbox{\,$\inbar\kern-.3em{\rm Q}$}}
\def\IF{\rlx\hbox{\rm I\kern-.18em F}}
\def\IG{\rlx\hbox{\,$\inbar\kern-.3em{\rm G}$}}
\def\IH{\rlx\hbox{\rm I\kern-.18em H}}
\def\II{\rlx\hbox{\rm I\kern-.18em I}}
\def\IK{\rlx\hbox{\rm I\kern-.18em K}}
\def\IL{\rlx\hbox{\rm I\kern-.18em L}}
\def\one{\hbox{{1}\kern-.25em\hbox{l}}}
\def\0#1{\relax\ifmmode\mathaccent"7017{#1}%
B        \else\accent23#1\relax\fi}

\begin{document}

\begin{titlepage}

\vspace{.2in}
\begin{center}
{\large\bf Some results on natural numbers represented by quadratic polynomials in two variables}
\end{center}

\vspace{.2in}

\begin{center}

B. M. Cerna Magui\~na$^{(a)}$, H. Blas$^{(b)}$ and V. H. L\'opez Sol\'{\i}s$^{(a)}$

\par \vskip .2in \noindent

$^{(a)}$ Carrera Profesional de Matem\'atica\\
Universidad Nacional Santiago Ant\'unez de Mayolo\\
Campus Shancay\'an, Av. Centenario 200, Huaraz - Per\'u.\\
$^{(b)}$Instituto de F\'{\i}sica\\
Universidade Federal de Mato Grosso\\
Av. Fernando Correa, $N^{0}$ \, 2367\\
Bairro Boa Esperan\c ca, Cep 78060-900, Cuiab\'a - MT - Brazil. \\
  
\normalsize
\end{center}

\vspace{.3in}

\begin{abstract}
\vspace{.3in}
We consider a set of equations of the form $p_j (x,y) = (10 x+m_j)(10 y + n_j),\,\,x\geq 0, y\geq0$, $j=1,2,3$, such that   $\{m_1=7, n_1=3\}$, $\{m_2=n_2=9\}$ and $\{m_3=n_3=1\}$, respectively. It is shown that if $(a(p_j),b(p_j)) \in \IN \times \IN$ is a solution of the $j'$th equation one has the inequality $\frac{p_j}{100}\leq A(p_j) B(p_j) \leq \frac{121}{10^4} p_j$, where  $A(p_j)\equiv a(p_j)+1, B(p_j)\equiv b(p_j)+1\,$ and $p_{j}$ is a natural number ending in 1,  such that $\{A(p_1)\geq 4, B(p_1)\geq 8\}$, $\{A(p_2) \geq 2,  B(p_2)\geq 2\}$, and  $\{A(p_3) \geq 10, B(p_3)\geq 10\}$ hold, respectively.  Moreover, assuming the previous result we show that $1\leq \( \frac{A(p_j+10) B(p_j+10)}{A(p_j) B(p_j)}\)^{1/100} \leq e^{0,000201} \times (1+ \frac{10}{p_j})^{(0,101)^2}$, with  $\{A(p_1)\geq 31, B(p_1)\geq 71\}$, $\{A(p_2) \geq 11,  B(p_2)\geq 11\}$, and  $\{A(p_3) \geq 91, B(p_3)\geq 91\}$, respectively. Finally, we present upper and lower bounds for the relevant positive integer solution of the equation defined by $p_j = (10 A+m_j)(10 B + n_j)$, for each case $j=1,2,3$, respectively.
\end{abstract}

{\parindent= 4em \small  \sl Keywords: Primes, Diophantine equation, natural numbers, quadratic polynomials.}

\end{titlepage}

\section{Introduction}

Number theory, and in particular the theory of  prime numbers, still fascinates mathematicians and recently the physicists (see e.g. \cite{r1, r31} and references therein). 
Nowadays, the data sets come from computer algorithms (see e.g. \cite{r11}  for a highly optimized sieve of Eratosthenes implementation, which counts the primes below $10^{10}$ in just 0.45 seconds), but mathematicians are still pursuing new patterns in primes.  Recently, it has been shown that, contrary to the commom believe that the primes occur randomly, prime numbers have a peculiar dislike for other would-be primes which end in the same digit \cite{r3}. In fact, when performed a randomness check on the first 100 million primes, they found that a prime ending in 1 was followed by another prime ending in 1 only 18.5 percent of the time - different from the 25 percent you would expect given that primes greater than five can only end in one of four digits: 1, 3, 7, or 9. The current research focuses on finding new patterns in prime numbers (see e.g. \cite{r3, r31} and references there in), as well as primes with largest number of digits (as of July 2018, the largest known prime is $2^{77,232,917}  - 1$ with nearly 22 million digits long.  It was found by the Great Internet Mersenne Prime Search (GIMPS) \cite{r33}). 

Fermat-Euler prime number theorem states that a prime number $p > 2$ is a sum of two squares of
integers, if and only if $p \equiv 1 (\mbox{mod}\, 4)$. In particular, $m^2 + n^2 \in \IP$ for infinitely
many integers $m, n$. H. Iwaniec generalized \cite{iwaniec} to 
polynomials of degree $2$, in two variables, satisfying certain assumptions. Let $P(m, n) = a\, m^2 + b\, m n + c\, n^2 + e\, m + f\, n + g$
be a primitive polynomial with integer cofficients. For $P(m, n)$ reducible
in $\IQ[m, n]$ the Dirichlet's theorem on arithmetic pogression can be used to answer the question whether it represents infinitely many primes \cite{luca}. If $P(m, n)$ is
irreducible the theorem in \cite{iwaniec} can be applied.
In our previous contribution \cite{r2}, we have discussed some aspects of the problem of obtaining a prime number starting from a given prime number. 

In this paper, following  our previous result on prime numbers we provide two new results, which together with the Theorem 2 of \cite {r2} would be relevant for constructing primes with millions of digits. The first result establishes that if $(a(p),b(p)) \in \IN\times \IN$ is a solution of any of the equations $i)\,\,  p = (10 x +7 ) (10 y + 3);
ii) \,\,  p = (10 x +9 ) (10 y + 9); iii)  \,\, p = (10 x +1) (10 y + 1),\,\, x\geq 0,\,y\geq 0$ then $\frac{p}{100} \leq A B \leq \frac{121}{10^4} p$, where $A(p)=a(p)+1,\,B(p)=b(p)+1$ and  $p$ is a natural number ending in 1, such that for the case $i)  \,\, A  \geq 4, \,\,\,\,\, B \geq 8$, case
$ii)  \,\, A \geq   2, \, \,\,\,\, B \geq 2$ and case $iii)  \,\, A \geq 10, \,\,\, B \geq 10$.  

The second result assumes the above result as a hypothesis and establishes that $1\leq \( \frac{A(p+10) B(p+10)}{A(p) B(p)}\)^{1/100} \leq e^{0,000201} \times (1+ \frac{10}{p})^{(0,101)^2}$, where  $A(p)=a(p)+1,\,B(p)=b(p)+1$, such that for the case $i)  \,\, A  \geq 31, \,\,\,\,\, B \geq 71$, case
$ii)  \,\, A \geq   1, \, \,\,\,\, B \geq 1$ and case $iii)  \,\, A \geq 91, \,\,\, B \geq 91$.

Moreover, we establish upper and lower bounds for the positive integer solutions of the equation $ p = (10 x +7 ) (10 y + 3)$. Analogous results are obtained for the remaning cases $ii)$ and $iii)$.

In this work $\IN$ represents the set of natural numbers, $e$ is the Euler's number, $\IN\times\IN = \{ (a,b)/ a \in \IN, b\in \IN \}$, $[X(p)]'$ stands for $p-$derivative of $X(p)$, and $\frac{1}{1-x}=\sum_{k=1}^{\infty} x^k, \, |x|< 1$.
  
Also, denote by $a(p)$ and $b(p)$ natural numbers that depend on $p$, $p\in \IN$.

\section{Some natural numbers ending in $1$ and quadratic polynomials.}

{\bf Theorem 1 }.  Let $p$ be a fixed natural number ending in 1.  If there exist $(a,b) \in \IN \times \IN$ such that one of the following relationships hold
\br
i)\,\,  p &=& (10 a +7 ) (10 b + 3) \notag\\
ii) \,\,  p &=& (10 a +9 ) (10 b + 9) \notag\\
iii)  \,\, p &=& (10 a +1) (10 b + 1), \notag
\er
then, one has 
\br
\frac{p}{100} \leq A B \leq \frac{121}{10^4} p, \label{th0}
 \er
where $A=a+1$, $B=b+1$, with the next conditions for the relevant cases above
\br
i)  \,\, A & \geq &4, \,\,\,\,\, B \geq 8  \notag\\
ii)  \,\, A &\geq  & 2, \, \,\,\,\, B \geq 2  \notag\\
iii)  \,\, A &\geq &10, \,\,\, B \geq 10. \notag
\er
{\bf Demonstration.}

Let us consider the case $ i) $.  So take $p = (10 x +7 ) (10 y + 3)$,  with  $x\geq 0,\,y\geq 0$. Next we take $p = (10 X -3 ) (10 Y -7)$, where $X=x=1,\, Y=y+1 $. Thus, one has
\br
p&=&100 X Y (1-\frac{3}{10 X})(1- \frac{7}{10 Y}),\,\,\, X \geq 1, Y \geq 1, \label{e1}
\er
which implies
\br
p &\leq& 100 X Y.\label{th11}
\er
Furthermore, by $(\ref{e1})$
\br
p \frac{1}{1-\frac{3}{10 X}} \frac{1}{1- \frac{7}{10 Y}} &=& 100 X Y.\er \label{e2}
Since $0 < \frac{3}{10 X} < 1,\,\, 0 < \frac{7}{10 Y}< 1$ the equation $(2.4)$ can be written as
\br
100 X Y &=& p (1+ \frac{3}{10 X}+ \frac{9}{100 X^2}+...)(1+\frac{7}{10 Y}+ \frac{49}{100 Y^2}+...) \label{th2}
\er
One can write the next identities for the series expansions above
\br
\frac{3}{10 X}+ \frac{9}{100 X^2}+\frac{27}{1000 X^3}+...&=& \frac{3}{10 X-3} \leq \frac{1}{10}, \,\,\,{\bf for}\,\, X > \frac{33}{10},  \label{e3}\\
\frac{7}{10 Y}+ \frac{49}{100 Y^2}+\frac{343}{1000 Y^2}+... &=&\frac{7}{10 Y-7} \leq \frac{1}{10},\,\,\,{\bf for}\,\, Y > \frac{77}{10}.   \label{e4}
\er
Taking into account these relationships $(\ref{e3})$ and $(\ref{e4})$ one has that  the equation (\ref{th2}) can be written as
\br
100 X Y \leq \frac{121}{100} p \label{e5}.
\er
Therefore, from $(\ref{th11})$ and  $(\ref{e5})$ one has
\br
 \frac{p}{100} \leq X Y \leq \frac{121}{10^4} p. \notag
\er
So, for natural numbers $X, Y$ one  has shown that the equation (\ref{th0}) holds.

{\bf Observation 1.} The demonstrations for the cases $ii)$ and $iii)$ follow analogous steps.
 
{\bf Observation 2.} Notice that  another relationship for the cases $i), ii)$ and $iii)$ becomes
\br
100 a b \leq p \leq 100 (a+1) (b+1). \notag
\er
 
{\bf Observation 3.} The positive integer number solutions of the Diophantine equation $p=(10 x + 7) (10 y + 3)$, where $p$ is a natural number ending in 1, take the following forms
\br
(a,b) & =& (\frac{p-21-10 \lambda}{30}, \frac{3 \lambda}{p-10 \lambda})   \notag\\
\mbox{or} &&    \notag\\
 (a,b)& = &(\frac{7 \lambda}{p-10 \lambda}, \frac{p-21-10 \lambda}{70}),\,\, \mbox{for some}\,\, \lambda \in \IN.  \notag
 \er

{\bf Theorem 2}.  Let $p$ a natural number ending in 1.  Assume that there exist $(a(p),b(p)) \in \IN \times \IN$ such that one of the following relationships hold
\br
i)\,\,  p &=& (10 a(p) +7 ) (10 b(p) + 3)   \notag\\
ii) \,\,  p &=& (10 a(p) +9 ) (10 b(p) + 9) \notag \\
iii)  \,\, p &=& (10 a(p) +1) (10 b(p) + 1).   \notag
\er
Then, one has 
\br
1\leq \( \frac{A(p+10) B(p+10)}{A(p) B(p)}\)^{1/100} \leq e^{0,000201} \times (1+ \frac{10}{p})^{(0,101)^2}, \label{e7}
\er
where $A(p)=a(p)+ 1,\, B(p)=b(p)+ 1$, with the following  condition for each case above
\br
i)  \,\, A(p) & \geq &31, \,\,\,\,\, B(p) \geq 71   \notag\\
ii)  \,\, A(p) &\geq  & 11, \, \,\,\,\, B(p) \geq 11    \notag\\
iii)  \,\, A(p) &\geq &91, \,\,\,\, B(p) \geq 91.   \notag
\er

{\bf Demonstration.} The product $X(p)Y(p)$ of the solutions $X(p)$ and $Y(p)$ in each case  $i)$, $ii)$ and $iii)$ can be arranged as a sequence of increasing numbers. Then, analogous to what was done in the demonstration of the Theorem $1 (i)$
\br
\frac{p}{100} \leq X(p) Y(p) \leq \frac{101^2}{10^6} p \label{th21}.
 \er
Since the product $X(p)Y(p)$ as an increasing function of $p$ then $[X(p)Y(p)]'$ is a positive number; so, taking into account  this observation one has that  (\ref{th21}) becomes 
\br
\frac{10^6}{101^2 p} [X(p)Y(p)]' \leq \frac{[X(p)Y(p)]'}{X(p)Y(p)} \leq \frac{100}{p} [X(p)Y(p)]'. \label{e6}
\er
Upon integrating the terms of the inequality $(\ref{e6})$ one has
\br
1\leq \( \frac{X(p+10) Y(p+10)}{X(p) Y(p)}\)^{1/100} \leq e^{\int_{p}^{p+10} \frac{[X(p)Y(p)]' }{p} dp} \label{th22}.
\er
We need some identities. Integrating by parts the exponential of the r.h.s. of the inequality $(\ref{th22})$ one has
\begin{equation}
\begin{split}
&I \equiv \int_{p}^{p+10} \frac{[X(p)Y(p)]' }{p} dp = \frac{X(p+10)Y(p+10)}{p+10} - \frac{X(p)Y(p)}{p}\\[4pt]
&~~~~~~~~~~~~~~~~~~~~~~~~~~~~~~~~~~~~~+\int_{p}^{p+10} \frac{X(p)Y(p)}{p^2} dp\label{th23}.
\end{split}
\end{equation}
Next, from (\ref{th21})  one can write
\br
\frac{1}{100} & \leq & \frac{X(p+10)Y(p+10)}{p+10} \leq \frac{101^2}{10^6}\label{th24}\\
-\frac{101^2}{10^6} & \leq &-\,\, \frac{X(p)Y(p)}{p} \,\,\,\,\,\,\,\,\,\,\,\leq -\frac{1}{100} \label{th25}.
\er
So from (\ref{th23}) and  (\ref{th24})-(\ref{th25}) we have 
\br
I \leq \frac{201}{10^6} + \int_{p}^{p+10} \frac{X(p)Y(p)}{p^2}  dp \label{th26}.
\er
Then, from (\ref{th21}) and the last relationship (\ref{th26})  
\br
I \leq 0,000201 + (0,101)^2 \ln(\frac{p+10}{p}). \label{th27}    
\er
Finally, from (\ref{th22}), (\ref{th23}) and the last relation (\ref{th27}) we obtain $(\ref{e7})$
\br
1\leq \( \frac{X(p+10) Y(p+10)}{X(p) Y(p)}\)^{1/100} \leq e^{0,000201} \times (1+ \frac{10}{p})^{(0,101)^2}. \notag
\er

{\bf Theorem 3}.  Let $(a,b) \in \IN \times \IN$ be a solution of the equation 
\br\label{th3}
p+10 = (10 x + 7)(10 y + 3),\er 
where $p$ is a prime number ending in the unity. If $(x_0 , b) \in \IQ^{+}\times \IN$ is a solution of the equation (\ref{th3}), then one has
\br
1 < \frac{A}{X_0} < (1+ \frac{10}{p}) \times \frac{101}{100},\,\,\, \,\, A\geq 31, \label{e10}
\er
\br
    \frac{10(p+10)}{11 p} < \frac{A}{X_0} <  \frac{101(p+10)}{100 p},\,\,\,\, \forall \, X_0 \geq 3.3, \label{e11}
\er
where $A = a+1,\,X_0=x_0 + 1$.

{\bf Demonstration.} Taking into account the hyphotesis of the Theorem we can write  from (\ref{th3}) the next identity
\br
\frac{p+10}{p} = \frac{10 A -3}{10 X_0 -3}. \label{e8}
\er 
Expanding in power series, the equation $(\ref{e8})$ becomes
\br
X_0 (p+10)(1+\frac{3}{10A}+\frac{9}{100 A^2}+...) = A p (1+\frac{3}{10X_0}+\frac{9}{100 X_0^2}+...).\label{e9}
\er
From the identity $(\ref{e9})$ one can get the inequalities $(\ref{e10})$ and $(\ref{e11})$
\br
1< \frac{A}{X_0}  < (1+ \frac{10}{p}) \times \frac{101}{100},\,\,\,\, \, \forall  A \geq 31 
\er
and 
\br
 \frac{10(p+10)}{11 p} < \frac{A}{X_0} <  \frac{101(p+10)}{100 p},\,\,\,\, \forall \, X_0 \geq 3.3.
\er
where $A = a+1,\,X_0=x_0 + 1$.

{\bf Theorem 4}.  Let $(a,b) \in \IN \times \IN$ be  integer solutions of the equation 
\br\label{th4}
p= (10 x + 7)(10 y + 3),  \,\, x\geq 0,\,\, y\geq 0\er 
Then
\br
1 \leq  \frac{p+70 A -21}{p+259} \leq e^{(7 \sqrt{p}-259)/p}\,\, \,\,\mbox{if}\, \,\, 4 \leq A\leq \frac{\sqrt{p}+3}{10}  \label{e13}
\er
or
\br
1 \leq  \frac{p+30 B -21}{p+219} \leq e^{(3 \sqrt{p}-3)/p}\,\, \,\,\mbox{if}\,\,\,  8 \leq B \leq \frac{\sqrt{p}+7}{10}, \label{e14}
\er
where $A = a+1,\, B=b + 1$.

{\bf Demonstration.} Taking into account the Theorem 1 one can write 
\br
\frac{p}{100} \leq X Y(X) \leq \frac{121}{10^4} p,\,\,\,\,X \geq 3.3,\,\,\,\,Y\geq 7.7. \label{e12}
\er
From the relation $(\ref{e12})$ and $p=(10 X -3)(10Y -7)$ one has
\br
\frac{10^3}{121p} X \leq  \frac{10X-3}{p+70X-21} \leq \frac{10}{p} X\label{th41}
\er
Upon integrating the last inequality in the interval $X=[4, A]$ one has
\begin{equation}
\begin{split}
& \exp{\{\frac{70(A-4)}{p}-\frac{2450}{p^2} (A^2-16)\}} \leq \frac{p+70A-21}{p+259}\\[4pt]
&~~~~~~~~~~~~~~~~~~~~~~~~~~~~~~~~~~~~~~~~~~~~~\leq \exp{\{\frac{70}{p} (A-4)-\frac{49 \times 10^4 (A^2-16)}{242}\}} \label{th42}.
\end{split}
\end{equation}
From (\ref{th42}) and the condition  $ 4 \leq A\leq \frac{\sqrt{p}+3}{10}$ we have $(\ref{e13})$ and $(\ref{e14})$
\br
1 \leq \frac{p+70A-21}{p+259} \leq e^{\frac{7 \sqrt{p}-259}{p}},\,\,\,\, \mbox{if} \,\,\,\,4 \leq A\leq \frac{\sqrt{p}+3}{10}.    \notag
\er  

Similarly, we can get
\br
1 \leq \frac{p+30B-21}{p+219} \leq e^{\frac{3 \sqrt{p}-3}{p}},\,\,\,\, \mbox{if} \,\,\,\,8 \leq B\leq \frac{\sqrt{p}+7}{10}.    \notag
\er

{\bf Corollary}.  Let $p$ be a natural number ending in $1$. So, we have the following three cases.

$i) \,\, \mbox{If}\,\, p = (10x+7)(10y+3)$ has positive integer solution then there are integers $(x,y)$ in the region $[3, + \infty> \times [7,+\infty>$  such that
\br
x&=&\frac{d-4}{14}\pm \frac{5^{k-1} 2^{j-1}}{7},   \notag\\
y&=&\frac{d+4}{6}\mp \frac{5^{k-1} 2^{j-1}}{3}, \,\,\,\,\, \mbox{for}\,\, d \,\,\mbox{even},  \mbox{and} \,\,\,k\geq 1,\, j\geq1,  \notag
\er
or
\br
x&=&\frac{d-4}{14}\pm \frac{5^{k-1} \tau^j}{14},  \notag \\
y&=&\frac{d+4}{6}\mp \frac{5^{k-1} \tau^{j}}{6}, \,\,\,\,\, \mbox{for}\,\, d\,\, \mbox{and}\, \tau \,\, \mbox{odd integers},\, \,\mbox{and} \,\,  k\geq   1,\, j\geq1,  \notag
\er 
where $d \equiv 7x + 3y$ with $\frac{\sqrt{21p}-29}{5} \leq d \leq 0.021p-7.9$.

$ ii) \,\, \mbox{If}\,\, p = (10x+9)(10y+9)$ has positive integer solution then there are integers $(x,y)$ in the region $[2, + \infty> \times [2,+\infty>$  such that
\br
x&=&\frac{d}{2}\pm 5^{k-1} 2^{j-1},\,   \notag\\
y&=&\frac{d}{2}\mp 5^{k-1} 2^{j-1},\, \, k,j \in \IN,\,\, k\geq 1,\, j\geq 1,  \notag
\er
or
\br
x&=&\frac{d}{2}\pm \frac{5^{k-1} \tau^j}{2},\,    \notag\\
y&=&\frac{d}{2}\mp \frac{5^{k-1} \tau^j}{2},\, \, k, j \in \IN,\,\, k\geq 1,\, j\geq 1, \mbox{and}\, \tau\, \mbox{odd}.   \notag
\er
where $d \equiv x + y$ with $\frac{\sqrt{p}-9}{5} \leq d \leq \frac{p-81}{90}$.

$ iii) \,\, \mbox{If}\,\, p = (10x+1)(10y+1)$ has positive integer solution then there are integers $(x,y)$ in the region $[9, + \infty> \times [9,+\infty>$ taking the same forms as the ones discussed in case $ii)$,  
where $\frac{\sqrt{p}-1}{5} \leq d \leq \frac{p-1}{10};\, d \equiv x+y, \,\, d\in \IN $.

{\bf Demonstration.} The proof in the case $i)$ follows by applying the Theorem 1. In fact, from (\ref{th0}) and   $ p = (10X(p)-3)(10Y(p)-7)$ one can get $0\leq 7x+ 3y \leq 0.021p - 7.9.$ As usual one takes $X(p)=x(p)+1,\,Y(p)=y(p)+1$.  

In addition, defining $d \equiv 7x + 3y$ and form the eq. $p = (10x+7)(10y+3)$ one can get $x= \frac{d-4}{14} \pm \frac{\sqrt{(5d+29)^2-21p}}{14\times 5}$. It is clear that in order to have positive integer solutions one must have 
\br
(5d+29)^2 =\left\{ \begin{array}{lcl} 21 p+5^{2k} 2^{2j}, & \mbox{for} & d \,\,\, \mbox{even},\, k,j \in \IN \\
21 p+5^{2k} \tau^{2j} ,& \mbox{for} & d\,\,\, \mbox{odd},\, k, j \in \IN,\, \mbox{and}\, \tau \, \mbox{odd}, \end{array}  \notag
\right.
\er 
where  $\frac{\sqrt{21p}-29}{5} \leq d \leq 0.021p-7.9$. 

Analogous steps can be followed for the proofs in the cases $ii)$ and $iii)$.

\section{Acknowledgments}

BMCM would like to thank Concytec for partial financial support. HB would like to thank the members of the FC-UNASAM for hospitality and FC-UNI for partial financial support. VHLS thanks FC-UNASAM for hospitality.

\end{document}